\documentclass[]{amsart}
\usepackage[foot]{amsaddr}
\usepackage[a4paper,top=2.5cm,bottom=3cm,inner=3cm,outer=3cm]{geometry}

\usepackage[utf8]{inputenc}

\usepackage{amsmath,amsthm}
\usepackage{amsfonts,amssymb,mathrsfs}

\usepackage{xy}
\usepackage[bookmarks=false,breaklinks=true]{hyperref}
\usepackage{float}
\usepackage{graphicx}
\usepackage[percent]{overpic}
\usepackage{color}

\usepackage{epigraph}
\definecolor{myurlcolor}{rgb}{0,0,0.7}
\usepackage{hyperref}
\hypersetup{colorlinks,
linkcolor=myurlcolor,
citecolor=myurlcolor,
urlcolor=myurlcolor}
\usepackage[mathscr]{euscript}

\input xy
\xyoption{all}

\newcommand{\maps}{\colon}    

\newcommand{\Z}{{\mathbb Z}}  
\newcommand{\R}{{\mathbb R}}  
\newcommand{\C}{{\mathbb C}}  
\renewcommand{\H}{{\mathbb H}}  

\newcommand{\End}{{\rm End}}
\newcommand{\GL}{{\rm{GL}}}     

\theoremstyle{definition}

        \newcommand{\be}{\begin{equation}}
        \newcommand{\ee}{\end{equation}}
        \newcommand{\ba}{\begin{eqnarray}}
        \newcommand{\ea}{\end{eqnarray}}
        \newcommand{\ban}{\begin{eqnarray*}}
        \newcommand{\ean}{\end{eqnarray*}}
        \newcommand{\barr}{\begin{array}}
        \newcommand{\earr}{\end{array}}

\begin{document}
\title{The Tenfold Way}
\author[Baez]{John C.\ Baez} 
\address{Department of Mathematics, University of California, Riverside CA, 92521, USA}
\address{Centre for Quantum Technologies, National University of Singapore, 117543, Singapore}
\date{\today}
\maketitle

\epigraph{
\emph{Three Rings for the Elven-kings under the sky, \\
Seven for the Dwarf-lords in their halls of stone....}}
{J.\ R.\ R.\ Tolkien}

The tenfold way burst into prominence around 2010 when it was applied to physics: 
in simple terms, it implies that there are ten fundamentally different kinds of matter \cite{RSFL}.  But the underlying mathematics is much older.   

The three real division algebras---the real numbers $\R$, the complex numbers $\C$, and the quaternions $\H$---show up naturally whenever groups
act as linear transformations of real or complex vector spaces.   This fact is especially important in quantum mechanics, which describes physical systems using linear algebra.  One consequence is that particles can be classified into three kinds.  The ramifications also percolate throughout mathematics.   In 1962, Freeman Dyson called this circle of ideas the ``threefold way'' \cite{Dyson}.

The story becomes even more interesting when we study $\Z/2$-graded vector spaces, also known as ``super vector spaces''.  A super vector space is simply a vector space written as a direct sum of two parts, called its even and odd part.  In physics, these parts can be used to describe two fundamentally different classes of particles: bosons and fermions, or alternatively, particles and antiparticles.  For the pure mathematician, super vector spaces are a fundamental part of the landscape of thought.  For example any chain complex, exterior algebra, or Clifford algebra is automatically a super vector space. 

In 1964, C.\ T.\ C.\ Wall \cite{Wall} studied real ``super division algebras'', which we will explain in a moment.  He found that besides $\R, \C$ and $\H$, which give ``purely even'' super division algebras, there are seven more!   And thus the tenfold way was born.

But we should begin at the beginning.  The story starts with Hamilton, who in 1843 was the first to deliberately create an algebra, namely the quaternions.    This is a four-dimensional real vector space
\[    \H = \{a1 + bi + cj + dk : \; a,b,c,d \in \R \} \]
with the associative bilinear product determined by these relations:
\[     i^2 = j^2 = k^2 = ijk = -1 \]
together with the requirement that $1$ be the multiplicative identity.

The quaternions are better than a mere algebra: they are a ``division algebra'', meaning one where every nonzero element has a multiplicative inverse.     As soon as Hamilton announced the quaternions, others rushed in and created more algebras, but in 1877 Frobenius proved that the only finite-dimensional division algebras over the field of real numbers are $\R, \C$ and $\H$.   The octonions come close, but they are not associative, and thus not part of today's story: we shall only consider associative algebras.  Also, every vector space from now on will be finite-dimensional.

Why are division algebras important?   One reason comes from group representation theory. 
Suppose you have a ``representation'' $\rho$ of a group $G$ on a vector space $V$: that is, a
homomorphism from $G$ to the group of invertible linear transformations of $V$.  
The linear transformations of $V$ that commute with every map $\rho(g) \maps V \to V$ 
form an algebra with composition as multiplication, called $\End(\rho)$.  Now suppose $\rho$ is ``irreducible'', meaning that the only subspaces preserved by all the maps $\rho(g) \maps V \to V$  are $\{0\}$ and $V$ itself.    Then Schur's lemma holds: any transformation $T \in \End(\rho)$ must have kernel and range equal to $\{0\}$ or $V$, so $T$ must be either invertible or zero. Thus, 
$\End(\rho)$ is a division algebra!

If $\rho$ is a \emph{real} irreducible representation of a group, we thus have three choices: 
$\R$, $\C$, and $\H$.    In quantum physics, states of a physical system are described by 
vectors in a Hilbert space, and symmetries act as unitary operators on this Hilbert space.
Elementary particles---systems that have no parts---are described by irreducible representations.   Thus, if quantum physics used \emph{real} Hilbert spaces, elementary particles would come in three kinds.

In fact, quantum physicists use \emph{complex} Hilbert spaces. The only finite-dimensional division algebra over an algebraically closed field is the field itself, and $\C$ is algebraically closed.  Thus, $\End(\rho)$ must be $\C$ for any irreducible representation $\rho$ of a group on a complex vector space.   But remarkably, there is still a threefold classification of irreducible representations of groups on complex vector spaces, coming from the three real division algebras \cite{Baez}.  So, elementary particles really do come in three kinds!

The whole story so far generalizes to super vector spaces.  As mentioned, a ``super vector space' is simply a vector space $V$ that is written as a direct sum of two subspaces:
\[     V = V_0 \oplus V_1 .\]
We think of $0$ and $1$ as elements of $\Z/2$.  We say an element of $V$ is ``even'' if it lies in $V_0$, ``odd'' if it lies in $V_1$, and ``homogeneous'' if it is either even or odd.     A map of super vector spaces $f \maps V \to W$ is a linear map that preserves the grading: $f(V_0) \subseteq W_0$ and $f(V_1) \subseteq W_1$.  The tensor product of super vector spaces is defined by
\[    (V \otimes W)_i = \bigoplus_{j+k = i} V_j \otimes W_k .\] 
A ``super algebra'' is an algebra $A$ that is equipped with the structure of a super
vector space and whose multiplication, regarded as a map $m \maps A \otimes A \to A$, is a
map of super vector spaces.   So, a super algebra is just an algebra $A = A_0 \oplus A_1$ obeying these rules:
\begin{itemize}
\item The product of two even elements is even.
\item The product of an even and an odd element is odd.
\item The product of two odd elements is even.
\end{itemize}
For example, take $A = \R[x]$, let $A_0$ consist of the polynomials that are even functions,
and let $A_1$ consist of those that are odd functions.

 A ``super division algebra'' is a superalgebra such that every nonzero \emph{homogeneous} element has an inverse.  Just as division algebras arise naturally from irreducible representations,  super division algebras arise naturally from irreducible representations in the super world.  We can generalize everything about Schur's lemma to this case.  We will not give the details here, but the upshot is that we can define a representation $\rho$ of a ``supergroup'' $G$ on a super vector space $V$, and then there is a super algebra $\End(\rho)$ of linear maps $T \maps V \to V$ that ``supercommute'' with every map $\rho(g) \maps V \to V$.   When $\rho$ is irreducible in the appropriate sense, $\End(\rho)$ must be a super division algebra.   

This makes it important to classify super division algebras.  It turns out that over the real numbers there are ten.  This is not hard to prove.  Suppose $A$ is a real super division algebra.  Then its even part $A_0$ is closed under multiplication, and the inverse of an even element must be even, so $A_0$ is a division algebra.   If $A_1 = \{0\}$, this gives three options:
\begin{itemize}
\item The super division algebra with $A_0 = \R$, $A_1 = \{0\}$.
\item The super division algebra with $A_0 = \C$, $A_1 = \{0\}$.
\item The super division algebra with $A_0 = \H$, $A_1 = \{0\}$.
\end{itemize}
In short, the three real division algebras can be seen as ``purely even'' super division algebras.

What if $A_1 \ne \{0\}$?  Then we can choose a nonzero element $e \in A_1$.  Since 
it is invertible, multiplication by this element sets up an isomorphism of vector spaces
$A_0 \cong A_1$.   In the case $A_0 = \R$, we 
can rescale $e$ by a real number to obtain either $e^2 = 1$ or $e^2 = -1$. 
This gives two options:
\begin{itemize}
\item The super division algebra with $A_0 = A_1 = \R$ and an odd element $e$
with $e^2 = 1$.
\item The super division algebra with $A_0 = A_1 = \R$ and an odd element $e$
with $e^2 = -1$.
\end{itemize}

In the case $A_0 = \C$, note that the map $a \mapsto eae^{-1}$ 
defines an automorphism of $A_0$, which must be either the identity or complex conjugation.   
If $eae^{-1} = a$ then $e$ commutes with all complex numbers, so we can rescale it by a 
complex number to obtain $e^2 = 1$.   If $eae^{-1} = \overline{a}$ then $ea =\overline{a}e$, so $e(e^2) = (e^2)e$ implies that $e^2$ is 
real.  Rescaling $e$ by a real number, we can obtain either $e^2 = 1$ or $e^2 = -1$.
So, we have three options:
\begin{itemize}
\item The super division algebra with $A_0 = A_1 = \C$ and an odd element $e$
with $ei = ie$ and $e^2 = 1$.
\item The super division algebra with $A_0 = A_1 = \C$ and an odd element $e$
with $ei = -ie$ and $e^2 = 1$.
\item The super division algebra with $A_0 = A_1 = \C$ and an odd element $e$
with $ei = -ie$ and $e^2 = -1$.
\end{itemize}

Finally, in the case $A_0 = \H$, the map $a \mapsto eae^{-1}$ defines an automorphism of $A_0$---which, it turns out, must be conjugation by some invertible quaternion $q$: $eae^{-1} = qaq^{-1}$.   But this means $q^{-1}e$ commutes with all of $A_0$, so replacing $e$ by $qe^{-1}$ we can assume our automorphism is the identity.    In other words, now $e$ commutes with all of $A_0$, so $e^2$ does as well.  Thus $e^2 \in \H$ is actually real.  Rescaling $e$ by a real number we can obtain either $e^2 = 1$ or $e^2 = -1$.   So, we have two options:
\begin{itemize}
\item The super division algebra with $A_0 = A_1 = \H$ and an odd element $e$
that commutes with everything in $A_0$ and has $e^2 = 1$.
\item
The super division algebra with $A_0 = A_1 = \H$ and an odd element $e$
that commutes with everything in $A_0$ and has $e^2 = -1$.
\end{itemize}
These are the ten real super division algebras!

The real fun starts here.  Wall showed that these ten algebras are all real or complex Clifford algebras.  The eight real ones represent all eight Morita equivalence classes of real Clifford
algebras.  The two complex ones do the same for the complex Clifford algebras.  The tenfold way thus unites real and complex Bott periodicity.   The ten real super division algebras also correspond naturally to the ten infinite families of compact symmetric spaces.  For the mathematician, the best way to learn about all these things is Freed and Moore's paper ``Twisted equivariant matter'' \cite{FM}.  The classification of real super division algebras was also explained by Deligne and Morgan \cite{DM}, and by Trimble \cite{Trimble}.

\end{document}